\documentclass[10pt, article]{amsart}

\usepackage{parskip}

\usepackage{geometry} 
\usepackage{ae} 
\usepackage[T1]{fontenc}
\usepackage[cp1250]{inputenc}
\usepackage{amsmath}
\usepackage{amssymb, amsfonts,amscd,verbatim}
\usepackage{MnSymbol}
\usepackage{mathrsfs} 
\usepackage{xypic}

\usepackage[normalem]{ulem}
\usepackage{hyperref}
\usepackage{indentfirst}
\usepackage{latexsym}
\input xy
\xyoption{all}

\usepackage{amsmath}    

\usepackage{tikz}

\theoremstyle{plain}
\newtheorem{Pocz}{Poczatek}[section]
\newtheorem{Proposition}[Pocz]{Proposition}

\newtheorem{Theorem}[Pocz]{Theorem}
\newtheorem{Corollary}[Pocz]{Corollary}

\newtheorem{Lemma}[Pocz]{Lemma}

\newtheorem{Example}[Pocz]{Example}

\theoremstyle{definition}
\newtheorem{Definition}[Pocz]{Definition}

\theoremstyle{remark}
\newtheorem{Remark}[Pocz]{Remark}

\newcommand{\floor}[1]{\left \lfloor #1 \right \rfloor}

\numberwithin{equation}{section}

\author{Thomas ~ Weighill}
\address{Tufts University, Medford MA, USA}
\email{thomas.weighill@gmail.com}

\title{Lifting coarse homotopies}

\date{ \today
}


\begin{document}

\begin{abstract}
Coarse geometry, and in particular coarse homotopy theory, has proven to be a powerful tool for approaching problems in geometric group theory and higher index theory. In this paper, we continue to develop theory in this area by proving a Coarse Lifting Lemma with respect to a certain class of bornologous surjective maps. This class is wide enough to include quotients by coarsely discontinuous group actions, which allows us to obtain results concerning the coarse fundamental group of quotients which are analogous to classical topological results for the fundamental group. As an application, we compute the fundamental group of metric cones over negatively curved compact Riemannian manifolds.
\end{abstract}

\maketitle

\section{Introduction}
Coarse geometry studies those properties of spaces which are invariant at large scale -- under quasi-isometry or, more generally, coarse equivalence. The notion of coarse equivalence arises naturally in geometric group theory, for example, because the word-metric on a finitely generated group does not depend -- up to coarse equivalence -- on the choice of generating set. Moreover, if a group acts geometrically on a proper geodesic metric space, then it is coarsely equivalent to that space (this is the Milnor-S\v{v}arc Lemma). This allows one to talk freely about large-scale properties such as the number of ends or Gromov hyperbolicity as properties of both the group (independent of presentation) and any space it acts on geometrically. Another consequence of the Milnor-S\v{v}arc Lemma is that the fundamental group of a closed Riemannian manifold $M$ is coarsely equivalent to the universal cover $\tilde{M}$ of that manifold. Roe was led to study the coarse geometry of complete Riemannian manifolds (for example, universal covers of compact manifolds) in order to formulate a notion of ``coarse index'' \cite{RoeIndex}. This coarse index is a generalization of the index of a differential operator on a compact manifold (the subject of the celebrated Atiyah-Singer Index Theorem) to the case where the manifold is no longer compact. Roe's ideas led to the formulation of the coarse Baum-Connes Conjecture (see e.g. \cite{YuCoarseBC} for a statement), which was later shown to be false. However, when the conjecture holds for a finitely generated group $G$ then it implies the Novikov Conjecture for all compact manifolds $M$ with $\pi_1(M) = G$. Guoliang Yu proved in 1998 that the coarse Baum-Connes Conjecture holds for spaces of finite asymptotic dimension \cite{YuNovikov}, later improving his result to include all spaces which coarsely embed into Hilbert space \cite{YuEmbed}.

Many of the successes of coarse geometry have been the result of translating important ideas and techniques from topology to the world of coarse geometry. Asymptotic dimension, mentioned above, was introduced by Gromov in \cite{Gromov93} as a natural coarse version of Lebesgue covering dimension. Much of Roe's work in index theory revolves around coarse cohomology, a coarse version of Alexander-Spanier cohomology. Coarse versions of homotopy have also played an important role. Among the very first proofs of the coarse Baum-Connes Conjecture for certain spaces were a proof for manifolds with ``Lipschitz good covers'' by Guoliang Yu \cite{YuCoarseBC} and a proof by Higson and Roe for Gromov hyperbolic spaces \cite{HigsonRoeCoarseBC}, both of which used some notion of coarse homotopy equivalence. More recently, coarse homotopy theory is being developed by Bunke, Engel and others using the context of $\infty$-categories \cite{BunkeEngel}. 

In \cite{Dranishnikov00}, Dranishnikov proposed that it would be interesting to define the coarse fundamental group using the upper half-plane, but the coarse fundamental group was only recently formally defined in a paper by Mitchener-Norouzizadeh-Schick (hereafter MNS) \cite{MitchenerNorSchick} using ideas from Mitchener's coarse homology theory. In that paper, MNS also introduce a notion of coarse homotopy very similar to one used by Bunke and Engel in \cite{BunkeEngel}, a definition which we will use with only a small modification (Definition \ref{coarsedef}). Another large-scale version of the fundamental group is the fundamental group at infinity studied in geometric group theory (see e.g. Chapter 16 of \cite{Geoghegan}), which is central to the long-standing Semistability Conjecture for finitely presented groups. 

The goal of this paper is to contribute to the development of coarse homotopy theory with an eye towards applications in geometric group theory and higher index theory. In particular, we prove a coarse homotopy lifting property and build out some related theory around it. Recall that the the homotopy lifting property for covering spaces from topology states that given any covering map $\pi: E \to B$, any homotopy $f: X \times I \to B$ and any continuous map $f_0: X \to E$ such that the diagram of solid arrows below commutes, there is a homotopy $\tilde{f}: X \times I \to E$ making the diagram commute.
$$
\xymatrix{
X \ar[d]_{1_X \times 0} \ar[r]^{f_0} & E \ar[d]^{\pi} \\
X \times I \ar[r]_f \ar@{..>}[ru]^{\tilde{f}} & B \\
}
$$

Our main result (Theorem \ref{lifting}), which we will call the Coarse Lifting Lemma, replaces continuous maps with their coarse analogues, bornologous maps, and replacing $X\times I$ by a coarse cylinder object $I_p X$ in this statement. It states  that for a surjective bornologous map $\pi$ satisfying certain conditions and for any bornologous map $f_0$ and coarse homotopy $f$ such that the diagram of solid arrows below commutes, there is a coarse homotopy $\tilde{f}$ making the diagram commute. 
$$
\xymatrix{
X \ar[d]_{i_0} \ar[r]^{f_0}  & E \ar[d]^{\pi} \\
I_p X  \ar@{..>}[ru]^{\tilde{f}} \ar[r]_f   & B \\
}
$$

The class of surjective bornologous maps $\pi$ for which this result holds includes quotients by group actions satisfying certain conditions  (see Lemma \ref{XGscattered}). This allows us to mimic the application of covering space theory to fundamental groups found in topology (see for example Chapters 9 and 13 of \cite{Munkres}) in the coarse setting. In particular, we obtain a short exact sequence 
$$\xymatrix{
0 \ar[r] & \pi_{1,q}^{\mathrm{coarse}}(X) \ar[r] &
 \pi_1^{\mathrm{coarse}}(X/G) \ar[r] & G \ar[r] &  0.
}
$$
for certain kinds of group actions (Theorem \ref{XGcoarse1}). In this short exact sequence, $\pi_{1,q}^{\mathrm{coarse}}(X)$ is a kind of coarse fundamental group of $X$ defined using only those maps which are coarse when post-composed with the quotient map $q: X \to X/G$, and $\pi_1^{\mathrm{coarse}}(X/G) $ is the coarse fundamental group of $X/G$. 

A useful way to construct spaces with interesting large scale behaviour is to take a compact space $M$ and construct the metric cone $\mathcal{O}M$ over it. In Section \ref{cones} we show how to apply Theorem \ref{XGcoarse1} to obtain an isomorphism 
$$
\pi_{1}^\mathsf{coarse}(\mathcal{O}M) \cong \pi_1(M),
$$
between the coarse fundamental group of $\mathcal{O}M$ and the (topological) fundamental group of $M$ when $M$ is a compact Riemannian manifold of non-positive sectional curvature (Theorem \ref{main}). This result, while not surprising given the corresponding result for simplicial complexes in \cite{MitchenerNorSchick}, appears to be new. Though, as we note at the end of the paper, it may be possible to derive this result from the simplicial complex result using triangulations of manifolds. We hope to apply Theorem \ref{XGcoarse1} to a wider class of metric cones as well as more general spaces in the near future.

\section{Coarse homotopies}
We will work in the setting of metric spaces rather than the abstract setting of coarse spaces introduced by Roe \cite{Roe} (see also the large scale spaces of Dydak-Hoffland \cite{DH}), but all the definitions and results in this paper can be easily generalized to that setting. Since a (classical) homotopy is a continuous map $X\times I \to Y$, it is natural to define a coarse homotopy via a cylinder object (to replace $X \times I$) and a large-scale notion of continuous map. The latter notion is usually called a bornologous map.

\begin{Definition}
A map $f$ from a metric space $X$ to a metric space $Y$ is called \textbf{bornologous} if there is some function $\rho$ (which we will call a \textbf{control function} for $f$) such that for any $x, x' \in X$,
$$
d(f(x), f(x')) \leq \rho(d(x,x')).
$$
A map is called \textbf{coarse} if it is bornologous and the inverse image of a bounded set is a bounded set.
\end{Definition}

Many proofs in coarse geometry require some form of ``pasting lemma''. We prove one such pasting lemma now since it will be useful later.

\begin{Lemma}\label{geodesicpasting}
Let $X$ be a geodesic metric space with a decomposition $X = A_1 \cup A_2 \cup \cdots \cup A_N$ where each $A_i$ is a closed set. Let $f: X \to Y$ be a map to a metric space $Y$ such that each restriction $f_i: A_i \to Y$ is coarse with control function $\rho_i$. Then $f$ is coarse. 
\end{Lemma} 
\begin{proof}
Let $x,x'$ be points in $X$ with $d(x,x') \leq R$ and let $\gamma$ be a geodesic from $x$ to $x'$. Since $\gamma$ is continuous and the $A_i$ are closed, we can choose a sequence $A'_1,\ A'_2,\ldots, A'_m$ of distinct $A_i$ and a sequence of real numbers $0 = t_0 \leq t_1 \leq \ldots \leq t_{m} = d(x,x')$ such that $\gamma(t_0) \in A'_1$, $\gamma(t_m) \in A'_m$ and for each $1 \leq k \leq m-1$, $\gamma(t_k) \in A'_{k} \cap A'_{k+1}$. Then
\begin{align*}
d(f(x), f(x')) & \leq \sum_{0 \leq k \leq m-1}  d(f(\gamma(t_k)), f(\gamma(t_{k+1})))\\
& \leq \sum_{0 \leq k \leq m-1} S(t_{k+1} - t_{k}) \\
& \leq N \cdot S(R)
\end{align*}
where $S(t) = \max_i \rho_i(t)$. It follows that $N\cdot S(\cdot)$ is a control function for $f$.
\end{proof}

We will use the definition of coarse homotopy introduced by MNS in \cite{MitchenerNorSchick} with a slight modification (see Remark \ref{MNSdef}). Denote the space $[0,\infty) \subseteq \mathbb{R}$ by $\mathbb{R}_+$. 

\begin{Definition}
Let $X$ be a metric space and $p: X \to [-1, \infty)$ be a function. Then the \textbf{$p$-cylinder} is defined as 
$$
I_pX = \{(x,t) \in X \times  \mathbb{R}_+ \mid t \leq p(x) + 1\}.
$$
with the $\ell^2$ metric. 
\end{Definition}

In other words, a coarse cylinder is a cylinder which ``opens out'' in a controlled way as you go to infinity. A canonical example we will need later is $I_p( \mathbb{R}_+ )$ where $p$ is the map $x \mapsto x-1$, which is just the subset of the plane $\{(x,t) \mid t \leq x\}$. We will call this space the \textbf{metric cone over $[0,1]$} and denote it by $c([0,1])$ (see Figure \ref{c01}).

\begin{figure}
\begin{tikzpicture}
\draw[fill=gray!60!white, draw=none] (0,0)--(2,2)--(2,0); 
\draw[<->, thick] (-1, 0)--(3,0);
\draw[<->, thick] (0, -1)--(0,3);
\draw (0,0)--(3,3);
\fill[black]  (2.2, 1.1) circle (1pt);
\fill[black]  (2.7, 1.35) circle (1pt);
\fill[black]  (3.2, 1.6) circle (1pt);
\end{tikzpicture}
\caption{$c([0,1])$}
\label{c01}
\end{figure}

When the map $p$ is coarse (as it will be when defining coarse homotopies), there are the evident inclusions $i_0: X \to I_pX$ and $i_1: X \to I_pX$ which send $x$ to $(x,0)$ and $(x, p(x) + 1)$ respectively. In the case of $c([0,1])$, these are the maps $i_0(x) = (x,0)$ and $i_1(x) = (x,x)$.

\begin{Definition} \label{coarsedef}
Let $X$ and $Y$ be metric spaces. Then a \textbf{coarse homotopy} is a coarse map $H: I_pX \to Y$ for some coarse map $p: X \to [-1, \infty)$. We say that the coarse homotopy $H$ is \textbf{from $f$ to $g$} if $H\circ i_0 = f$ and $H  \circ i_1 =g$.
\end{Definition}

To get an intuitive picture of this definition, notice that for every $x \in X$, $p(x, \cdot)$ gives a coarse map from $[0, p(x)+1]$ to $Y$, which we can think of as a coarse path associated to $x$. Just as for a classical homotopy, the definition requires the collection of all such paths to fit together in some sense. 

\begin{Remark}\label{MNSdef}
The original definition of coarse homotopy in \cite{MitchenerNorSchick} requires the map $p$ in $I_pX$ to have codomain $\mathbb{R}_+$ when defining a $p$-cylinder. Unfortunately, this means that the ``coarse loops'' used in the coarse fundamental group (see Definition \ref{coarsefundamentalgroup}) are not considered coarse homotopies. In topology, paths are the same as homotopies between maps from the singleton space, and this allows the use of lifting results for lifting both homotopies and loops. We would like to do the same in this paper, and so in order for our Coarse Lifting Lemma to apply to coarse loops, we use a more general definition of coarse homotopy which allows $p$ to have codomain $[-1,\infty)$. Clearly any coarse homotopy from $I_pX$ where $p$ takes values only in $\mathbb{R}_+$ (i.e.~in the sense of MNS) is also a coarse homotopy in our sense. Conversely, any coarse homotopy $H: I_p X \to Y$ from $f$ to $g$ with $p: X \to [-1,\infty)$ gives rise naturally to a coarse homotopy $H': I_{p'}X \to Y$ in the sense of MNS from $f$ to $g$ where $p'(x) = p(x)+1$ and
$$
H'(x,t) = \begin{cases}
H(x,t) & t \leq p(x) + 1\\
g(x) & p(x)+1  < t \leq p'(x) + 1
\end{cases}
$$
This is just the concatenation of two homotopies, so it is coarse by Proposition \ref{concat}. Therefore our notion of two maps being coarsely homotopic coincides exactly with that of MNS, so in particular our definition of $\pi^1_\mathsf{coarse}(X)$ is equivalent to that of MNS. 
\end{Remark}

\begin{Example}
Recall that two coarse maps $f,g: X \to Y$ are called \textbf{close} if $\sup_x d(f(x), g(x)) < \infty$. It is easy to check that \emph{any} two close maps from a metric space are coarsely homotopic. 
\end{Example}

Paralleling the classical case, we say that a coarse homotopy $H: I_pX \to Y$ is a \textbf{coarse homotopy relative to $A \subseteq X$} if $H(a, t) = H(a,0)$ for all $a \in A$ and all $0\leq t \leq p(a) +1$. An important property of coarse homotopies is that they can be concatenated. The following result is a slight generalization of Lemma 2.5 in \cite{MitchenerNorSchick}.

\begin{Proposition}\label{concat}
Let $H: I_pX \to Y$ and $H': I_qX \to Y$ be two coarse homotopies such that the maps $H\circ i_1$ and $H' \circ i_0$ are close. Then the concatenation $H\ast H': I_{p+q+1} X \to Y$ defined by
$$
H\ast H' (x,t) = \begin{cases}
H(x,t) & t \leq p(x) + 1 \\
H'(x, t-p(x)-1) & p(x) + 1< t \leq p(x) + q(x) + 2 \\ 
\end{cases}
$$
is a coarse homotopy.
\end{Proposition}
\begin{proof}
Define $H''$ to concide with $H$ except with $H''(x,p(x)+1) = H'(x,0)$ for all $x$. Then $H''$ is close to $H$, so $H''$ is coarse. Now that $H'' \circ i_0 = H' \circ i_1$, it follows from the proof of Theorem 2.4 in \cite{MitchenerNorSchick} that $H'' \ast H$ is coarse (in fact, the proof is just for $p$ and $q$ with codomain $\mathbb{R}_+$ but the proof goes through for $p$ and $q$ with codomain $[-1, \infty)$ too). Finally, $H \ast H'$ is close to $H'' \ast H'$ and thus is also coarse as required. 
\end{proof}

\begin{Remark}
There is also another very similar definition of cylinder object and coarse homotopy in the literature, namely that introduced by Bunke and Engel in \cite{BunkeEngel}. For a space $X$ and two coarse maps $p_+: X \to [0, \infty)$ and $p_-: X \to (-\infty, 0]$, the authors define the coarse cylinder to be the space
$$
I_{p_-, p_+} X = \{(x,t) \in X \times \mathbb{R} \mid p_-(x) \leq t \leq p_+(x) \}
$$
with the $\ell^2$ metric (or any metric coarsely equivalent to it). There are inclusions $i_0: x \mapsto (x,p_-(x))$ and $i_1: x \mapsto (x, p_+(x))$ and so we can define coarse homotopies from $f:X\to Y$ to $g: X \to Y$ in the obvious way. Any Bunke-Engel coarse homotopy $H: I_{p_-, p_+} X \to Y$ from $f$ to $g$ gives rise to a coarse homotopy $H': I_{p}X \to Y$ from $f$ to $g$ in the sense of this paper where $p(x) = - p_-(x) + p_+(x) - 1$ and $H'(x,t) = H(x+p_-(x), t)$. This may not be a coarse homotopy in the sense of MNS, but we have already seen that one can be induced by $H'$ (Remark \ref{MNSdef}). Conversely, any coarse homotopy $H: I_pX \to Y$ in the sense of this paper (or MNS) gives rise to a Bunke-Engel coarse homotopy $H':  I_{p_-, p_+} X \to Y$ where $p_-(x) = -p(x)/2$ and $p_+(x) = 1+p(x)/2$ and $H'(x,t) = H(x+p(x)/2,t)$. So for all the definitions of coarse homotopy considered so far (Bunke-Engel, MNS and this paper's slight modification of MNS), the statements ``$f$ is coarsely homotopic to $g$'' are equivalent.
\end{Remark}

\section{Coarse Lifting lemma}
Now that we have a definition of coarse homotopy, we proceed to state and prove the Coarse Lifting Lemma. The classical Lifting Lemma is stated for covering maps, so we introduce a coarse analog here. 

\begin{Definition}
Let $f: X \to Y$ be a bornologous map between metric spaces. Then $f$ is called a \textbf{soft quotient map} if it is surjective and for every $R > 0$ there is an $S > 0$ such that if $d(f(x), y) \leq R$ for some $x \in X$, $y \in Y$, then there is an $x' \in f^{-1}(y)$ with $d(x,x') \leq S$. 
\end{Definition}

The terminology is based on \cite{DikranjanZava}, in which the authors define (weakly) soft maps in the context of balleans. The following example is also the main example for the applications in this paper.

\begin{Example}\label{exgroup}
Let $G$ be a group acting on a metric space $X$ by isometries, and let $X/G$ be the orbit space equipped with the metric 
$$
d([x], [y]) = \inf\{d(x, y') \mid y' \in [y] \}.
$$
Then the quotient map $q: X \to X/G$ is a soft quotient map.
\end{Example}

We now introduce the following new definition.

\begin{Definition}
Let $f: X \to Y$ be a map between metric spaces. We say that $f$ has \textbf{scattered fibres} if for every $R > 0$, there is a bounded set $K$ in $Y$ such that if $f(x) = f(x') \notin K$ for $x, x' \in X$ and $d(x,x') \leq R$, then $x = x'$.
\end{Definition}

Recall that a surjective map $f:X \to Y$ between metric spaces was called \textbf{asymptotically faithful}  in \cite{WillettYu} if for every $R > 0$ there is a bounded set $K \subset Y$ such that $f$ is an isometry on every $R$-ball not intersecting $f^{-1}(K)$. It is easy to check that any such map has scattered fibres, though it need not be a soft quotient map. We are now ready to state and prove the Coarse Lifting Lemma.

\begin{Theorem}[Coarse Lifting Lemma]\label{lifting}
Let $\pi: A \to B$ be a soft quotient map with scattered fibres, let $f: I_pX \to B$ be a coarse homotopy and let $f_0: X \to A$ be a coarse map such that the diagram of solid arrows below commutes. Then there is a coarse homotopy $\tilde{f}: I_pX \to A$ making the diagram commute.
\begin{equation}\label{liftingdiagram}
\xymatrix{
X \ar[d]_{i_0} \ar[r]^{f_0}  & A \ar[d]^{\pi} \\
I_p X  \ar@{..>}[ru]^{\tilde{f}} \ar[r]_f   & B \\
}
\end{equation}
Moreover, if $\tilde{f}$ and $\tilde{f}'$ are two coarse homotopies making the diagram commute, then there is a bounded set $K$ in $I_pX$ such that $\tilde{f} |_{I_pX \setminus K} = \tilde{f}' |_{I_pX \setminus K}$.
\end{Theorem}
\begin{proof}
The main idea, based on the topological situation, is to lift each ``coarse path'' $f(x, \cdot)$ individually and prove that this gives the right map. Let $\rho$ be a control function for $f$. Since $\pi$ is a soft quotient map, there is a $T > 0$ so that if $\pi(a)$ is within $\rho(1)$ of $\pi(a')$ then there is an $a'' \in \pi^{-1}(\pi(a'))$ such that $d(a, a'') \leq T$. In particular, consider the points $\pi(f_0(x)) = f(x,0)$ and $f(x,\varepsilon)$ for any $0 < \varepsilon \leq 1$. Since $d(f(x,0),f(x,\varepsilon)) \leq \rho(1)$, we can find a $\tilde{f}(x, \varepsilon)$  so that $\pi(\tilde{f}(x,\varepsilon)) = f(x, \varepsilon)$ and  $d(\tilde{f}(x, \varepsilon), f_0(x)) \leq T$. Continuing by induction and varying $x$ we can define $\tilde{f}$ in such a way that $\pi \circ \tilde{f} = f$ and $\tilde{f}(x,n)$ is within $T$ of $\tilde{f}(x,n+\varepsilon)$ whenever $n \in \mathbb{N}$ and $\varepsilon \leq 1$. This implies that $d(\tilde{f}(x, s), d(\tilde{f}(x,t)) \leq |\floor{s}-\floor{t}|\cdot T + 2\cdot T$.

Now let $\rho_0$ be the control function for $f_0$ and let $R > 0$. Since $\pi$ is a soft quotient map, there is an $S > \rho_0(R)$ such that if $d(\pi(a), b) \leq \rho_0(R)$ then there is an $a' \in \pi^{-1}(b)$ with $d(a, a') \leq S$. Denote this $S$ by $S(R)$ since it depends only on $R$. Since $\pi$ has scattered fibres, there is a bounded set $K$ in $B$ such that if $\pi(a)=\pi(a')$ and $a$ is within $2S+2T$ of $a'$, then $a = a'$. Because $f$ and $p$ are coarse we can choose a bounded set $L$ in $X$ such that $L \times \mathbb{R}_+ \cap I_p X$ contains $f^{-1}(K)$. We will show by induction that $\tilde{f}(x, t)$ and $\tilde{f}(x', t)$ are at most $S$ apart whenever $x, x' \notin L$, $d(x, x') \leq R$. Suppose that $t = n +\varepsilon$ for $\varepsilon \in (0,1]$ and $d(\tilde{f}(x,n), \tilde{f}(x',n)) \leq S$. Let $a$ be a point in $\pi^{-1}(f(x', t))$ which is at most $S$ away from $\tilde{f}(x, t)$. Then $a$ and $x'$ are at most $2T + S + d(\tilde{f}(x,n), \tilde{f}(x',n)) \leq 2T+2S$ apart, and thus $a = x'$ proving that $\tilde{f}(x, t)$ and $\tilde{f}(x', t)$ are at most $S$ apart. Combining with the base case $d(\tilde{f}(x,0), \tilde{f}(x',0)) \leq \rho_0(R) \leq S$, we obtain our result for all $t$. Let $B(L, R)$ be the $R$-neighborhood of the set $L$. Its image under $\tilde{f}$ is clearly within bounded distance of $f_0(L\times \{0\})$, so it is bounded; denote its diameter by $L(R)$. Putting everything together, we have that $d((x,t),d(x',t')) \leq R$ implies
\begin{align*}
d(\tilde{f}(x,t), \tilde{f}(x',t')) & \leq  d(\tilde{f}(x,t), \tilde{f}(x,t'))+d(\tilde{f}(x,t'), \tilde{f}(x',t')) \\ 
& \leq (R+1)\cdot T + 2\cdot T + \max(S(R), L(R))
\end{align*}
using the fact that either $x,x' \in B(L, R)$ or both $x,x' \notin L$. It follows that $\tilde{f}$ is coarse as required. 

We now turn to uniqueness. Suppose $\tilde{f}$ and $\tilde{f}'$ are two distinct maps satisfying the requirements with control functions $\rho$ and $\rho'$ respectively. Using the fact that $\pi$ has scattered fibres, let $J$ be a bounded set such that if $d(a,a') \leq \rho(1) + \rho'(1)$, $\pi(a) = \pi(a')$ and $a \neq a'$, then $a, a' \in K$. Let $K = K' \times \mathbb{R}_+ \cap I_p X$ be a subset containing the bounded set $f^{-1}(J)$. We claim that $\tilde{f} |_{I_pX \setminus K} = \tilde{f}' |_{I_pX \setminus K}$. Suppose for contradiction that $\tilde{f}(x,t) \neq \tilde{f}'(x,t)$ and $x \notin K'$. Since $\tilde{f}(x,0) = \tilde{f}'(x ,0)$, we may choose $t$ so that $\tilde{f}(x,n) = \tilde{f}'(x,n)$ where $t = n + \varepsilon$ for $\varepsilon \in (0,1]$. Then $\tilde{f}(x,t)$ and $\tilde{f}'(x,t)$ are at most $\rho(1) + \rho'(1)$ apart, a contradiction to the choice of $K'$ because $\pi \tilde{f}(x,t) = \pi \tilde{f}(x',t') = f(x,t)$. This completes the proof.
\end{proof}

\begin{Remark}
One may ask whether some version of the Coarse Lifting Lemma holds for maps $\pi$ which are only coarsely surjective (i.e.~the image of $\pi$ is finite Hausdorff distance from the whole of $B$) and not necessarily surjective. As stated, Theorem \ref{lifting} cannot possibly hold for all such maps since the image of $f$ may not be contained in the image of $\pi$ in diagram (\ref{liftingdiagram}). On the other hand, one can easily show that given a commutative diagram of the form of diagram (\ref{liftingdiagram}), with $f: I_pX \to B$ a coarse homotopy, $f_0: X \to A$ a coarse map and $\pi$ a coarsely surjective map such that the restriction $\pi: A \to \pi(A)$ is a soft quotient map with scattered fibres, there is a coarse homotopy $\tilde{f}: I_pX \to A$ making the diagram commute \emph{up to closeness}.  In our main application (group actions), the maps will be surjective and commutativity only up to closeness would complicate the proofs, so we instead chose to focus on the more restrictive version of Theorem \ref{lifting} given above. 
\end{Remark}

\section{Coarse fundamental group}
Along with their definition of coarse homotopy, MNS also introduce the notion of coarse fundamental groups in \cite{MitchenerNorSchick}. Since the classical Lifting Lemma can be used to compute topological fundamental groups, it is no surprise that we have an analogous situation in the coarse setting. 

If $f,g : X \to Y$ are two coarse maps between metric spaces, then they are \textbf{coarsely homotopic} if there exists a coarse homotopy between them, that is, there is a coarse map $H: I_p X \to Y$ such that $H\circ i_0 = f$ and $H \circ i_1 = g$. Given a metric space $X$, an \textbf{$\mathbb{R}_+$-basepoint} is a coarse map $b: \mathbb{R}_+ \to X$.

\begin{Definition}\label{defc01}
  The metric cone over the interval, denoted by $c([0,1])$, is defined as $\{(x,t) \mid t \leq x\} \subseteq \mathbb{R}^2$ with the inherited $\ell^2$-metric. We denote the boundary $\{ (x,t) \in c([0,1]) \mid t = x \textsf{ or } t = 0 \}$ by $\partial c([0,1])$. 
 \end{Definition}

The metric cone $c([0,1])$ is clearly identical to the cylinder $I_p(\mathbb{R}_+)$, where $p(x) = x-1$, so there are natural inclusions $i_0: \mathbb{R}_+ \to c([0,1])$ and $i_1: \mathbb{R}_+ \to c([0,1])$. 

\begin{Remark}
It appears that the definition of $c([0,1])$ we give in Definition \ref{defc01}  differs from the representation of $c([0,1])$ in \cite{MitchenerNorSchick} (where $x \geq y$ is replaced by $x \leq y$), but ours seems more natural if one is inclined to view $c([0,1])$ also as a kind of $p$-cylinder, and the two definitions yield the same space with the same boundary up to isometry.
\end{Remark}

\begin{Definition}\label{coarsefundamentalgroup}
Let $X$ be a metric space with a chosen $\mathbb{R}_+$-basepoint $b: \mathbb{R}_+ \to X$. Then the \textbf{$1^\mathrm{st}$ coarse homotopy group} $\pi_{1}^{\mathrm{coarse}}(X, b)$ is the set of all relative coarse homotopy (relative to $\partial c ([0,1])$) classes of coarse maps $\alpha: c([0,1]) \to X$ such that $\alpha \circ i_0 = \alpha \circ i_1 = b$. 
\end{Definition}

MNS describe how to concatenate two coarse loops as follows. Given two coarse maps $\alpha, \beta: c([0,1]) \to X$ for which $\alpha \circ i_1 = \beta \circ i_0$, one can construct the concatenation $\alpha \ast \beta: c([0,1]) \to X$ via
$$
\alpha \ast \beta (x,t) = \begin{cases}
\alpha(x,2t) & t \leq x/2 \\
\beta(x, 2t-x) & x/2 < t \leq x \\ 
\end{cases}
$$
MNS prove that the operation $\ast$ gives rise to a well-defined group structure on $\pi_1^{\mathrm{coarse}}(X, b)$ given by $[\alpha]\cdot [\beta] = [\alpha \ast \beta]$ (Proposition 4.8 in \cite{MitchenerNorSchick}). In particular, $[\alpha \ast \beta]$ does not depend on the chosen representatives $\alpha$ and $\beta$.

We will be interested in results for general bornologous maps, not just for coarse maps. An obvious obstacle here is that a bornologous map $f: X \to Y$ does not in general induce a map on coarse fundamental groups since it may not send coarse homotopies to coarse homotopies. We will thus need a relative version of coarseness. This is just a special case of the notion of coarse map between bornological coarse spaces as introduced in \cite{BunkeEngel}, but we will not need the general theory here.

\begin{Definition}
Let $f: X \to Y$ be a bornologous map. A map $g: W \to X$ is called \textbf{$f$-coarse} if $g$ is bornologous and $f \circ g$ is coarse. 
\end{Definition}

Note that if $f$ is bornologous, then any $f$-coarse map is necessarily a coarse map. For a bornologous map $f: X \to Y$, there is then a corresponding notion of $f$-coarse homotopy between two maps $g_0, g_1: W \to X$, namely that there is a $f$-coarse map $H: I_p W \to X$ such that $H \circ i_0 = g_0$ and $H \circ i_1 = g_1$.

\begin{Definition}
Let $X$ be a metric space with a chosen $\mathbb{R}_+$-basepoint $b: \mathbb{R}_+ \to X$ and a bornologous map $f: X \to Y$. Then the \textbf{$1^\mathrm{st}$ $f$-coarse homotopy group} $\pi_{1,f}^{\mathrm{coarse}}(X, b)$ is the set of all relative $f$-coarse homotopy (relative to $\partial c ([0,1])$) classes of $f$-coarse maps $\alpha: c([0,1]) \to X$ such that $\alpha \circ i_0 = \alpha \circ i_1 = b$. 
\end{Definition}

It is easy to check that the same group operations work for the $f$-coarse homotopy group. Any bornologous map $f: X \to Y$ now induces a group homomorphism
$$
f_\ast: \pi_{1,f}^{\mathrm{coarse}}(X, b) \to \pi_{1}^{\mathrm{coarse}}(Y, f\circ b)
$$
for any $\mathbb{R}_+$-basepoint $b$ in $X$. Our next theorem is a coarse lifting correspondence similar to the lifting correspondence in topological covering space theory (see e.g. Chapter 9 of \cite{Munkres}), and the proof is based on the topological one.

\begin{Proposition}[Lifting Correspondence]\label{correspondence}
Let $\pi: X \to Y$ be a soft quotient map with scattered fibres. Let $b: \mathbb{R}_+ \to Y$ be an $\mathbb{R}_+$-basepoint in $Y$ and $b'$ a lift of $b$ to $X$. Suppose that for any other lift $b'': \mathbb{R}_+ \to X$ of $b$, there is a $\pi$-coarse homotopy $H$ from $b'$ to $b''$. Then there is a canonical (once $b'$ is chosen) bijection between the right cosets of $\pi_{\ast}(\pi_{1,\pi}^{\mathrm{coarse}}(X, b'))$ in $\pi_1^{\mathrm{coarse}}(Y, b)$ and equivalence classes of liftings $b'': \mathbb{R}_+ \to X$ of $b$ under the equivalence relation $b'' \sim b'''$ if 
$$
\{ t \in \mathbb{R} \mid b''(t) \neq b'''(t) \}
$$
is bounded.
\end{Proposition}
\begin{proof}
If $\alpha: c([0,1]) \to Y$ represents a class in $\pi_1^{\mathrm{coarse}}(Y, b)$, then we can lift it by Theorem \ref{lifting} to $\tilde{\alpha} : c([0,1]) \to X$ such that $\tilde{\alpha} \circ i_0 = b'$. Define $\Phi(\alpha)$ to be the map $\tilde{\alpha} \circ i_1$. We will show that $\Phi$ gives the required bijection. It is easy to show using the uniqueness part of Theorem \ref{lifting} that $[\Phi(\alpha)]_\sim$ is well-defined, and moreover that $[\alpha] = [\beta] \in \pi_1^{\mathrm{coarse}}(Y, b)$ implies $\Phi(\alpha) \sim \Phi(\beta)$. If $\Phi(\alpha) \sim \Phi(\beta)$, with $\tilde{\alpha}$ and $\tilde{\beta}$ lifts of $\alpha$ and $\beta$ respectively, then $\tilde{\alpha} \ast \tilde{\beta}^{\ast}$ (the concatenation of $\alpha$ with the reverse of $\tilde{\beta}$) represents an element of $\pi_{1,\pi}^{\mathrm{coarse}}(X, b')$ whose image under $\pi$ is $\alpha \ast \beta^\ast$. It follows that $\alpha$ and $\beta$ are in the same right coset of $\pi_{\ast}(\pi_{1,\pi}^{\mathrm{coarse}}(X, b'))$ in $\pi_1^{\mathrm{coarse}}(Y, b)$. An easy argument shows the converse, so that $\Phi$ descends to an injection from cosets of $\pi_\ast(\pi_{1,\pi}^{\mathrm{coarse}}(X, b'))$ in $\pi_1^{\mathrm{coarse}}(Y,  b)$ to $\sim$-equivalence classes of lifts. To show surjectivity, note that by assumption any lift $b''$ of $b$ is connected to $b'$ by a $\pi$-coarse homotopy $H: I_p \mathbb{R}_+ \to X$. We may adapt $b''$ so that $b''(0) = b'(0)$ without changing its $\sim$-class. An easy adaptation of Lemma 2.6 in  \cite{MitchenerNorSchick} now shows that $H$ gives rise to a map $H': c([0,1]) \to X$ with $H' \circ i_0 = b'$ and $H' \circ i_1 = b''$. Since $H'$ is a lift of the coarse homotopy $\pi \circ H$ from $b$ to $b$, we have that $\Phi(\pi \circ H) = [b'']$ as required.
\end{proof}

In the setting of the Proposition \ref{correspondence}, if $\pi$ is actually a coarse map then the $\pi$-coarse homotopy groups coincide with the coarse homotopy groups, so we get the following.

\begin{Corollary}\label{correspondencecor}
Let $\pi: X \to Y$ be a coarse soft quotient map with scattered fibres. Let $b: \mathbb{R}_+ \to Y$ be an $\mathbb{R}_+$-basepoint in $Y$ and $b'$ a lift of $b$ to $X$. Suppose that for any other lift $b'': \mathbb{R}_+ \to X$ of $b$, there is a coarse homotopy $H$ from $b'$ to $b''$. Then there is a canonical bijection between the cosets of $\pi_{\ast}(\pi_{1}^{\mathrm{coarse}}(X, b'))$ in $\pi_1^{\mathrm{coarse}}(Y, b)$ and $\sim$-equivalence classes of liftings $b': \mathbb{R}_+ \to X$ of $b$. In particular, if $\pi_{1}^{\mathrm{coarse}}(X, b')$ is trivial, then there is a canonical (once $b'$ has been chosen) bijection between $\pi_1^{\mathrm{coarse}}(Y, b)$ and $\sim$-equivalence classes of lifts of $b$.
\end{Corollary}

\section{Group actions}
We now turn our attention to quotients by group actions. For simplicity, we will consider only actions by isometries, but the results apply slightly more generally to uniformly bornologous actions. Given a group $G$ acting on a metric space $X$ by isometries, we will assume that the quotient space $X/G$ is given the metric as in Example \ref{exgroup}.

\begin{Definition}
Let $G$ be a group acting on a metric space $X$. We say that $G$ acts \textbf{uniformly coarsely discontinuously} if for every $R > 0$, there is a bounded set $K$ such that if $x \notin \cup_{g \in G} g(K)$ and $g$ is not the identity, then $d(x, g \cdot x) > R$.
\end{Definition}

One large class of examples of uniformly coarsely discontinuous actions is actions on metric cones induced by a geometric action on the underlying space (see Definition \ref{conedef} for a definition and Lemma \ref{discon} for the proof of this fact). Note that the term ``uniformly coarsely discontinuously'' is based on (but different from) the definition of ``coarsely discontinuously'' in \cite{WeighillHigginbotham}, which is more suited to the study of warped spaces than orbit spaces.

\begin{Lemma}\label{XGscattered}
Let $G$ be a group acting on a metric space $X$ by isometries. Let $X/G$ be the orbit space and $q: X \to X/G$ the natural quotient map. Suppose further that the action of $G$ is uniformly coarsely discontinuous. Then $q$ has scattered fibres.
\end{Lemma}
\begin{proof}
This is obvious from the definition, once one notices that $q^{-1}(q(K)) = \cup_{g \in G} g(K)$. 
\end{proof}

We are now ready to prove a result concerning the first coarse fundamental group of $X/G$.

\begin{Theorem}\label{XGcoarse1}
Let $G$ be a group acting on a metric space $X$ by isometries. Let $X/G$ be the orbit space and $q: X \to X/G$ the natural quotient map. Suppose further that the action of $G$ is uniformly coarsely discontinuous. Let $b: \mathbb{R}_+ \to X/G$ be an $\mathbb{R}_+$-basepoint such that for any two lifts $b', b'': \mathbb{R}_+ \to X$, there is a $q$-coarse homotopy $H$ from $b'$ to $b''$. Then for any lift $b'$ of $b$, we have a canonical short exact sequence
$$
\xymatrix{
0 \ar[r] & \pi_{1,q}^{\mathrm{coarse}}(X, b') \ar[r]^{q_\ast} &
 \pi_1^{\mathrm{coarse}}(X/G, b) \ar[r] & G \ar[r] &  0.
}
$$
\end{Theorem}
\begin{proof}
The map $q$ is a soft quotient map with scattered fibres, so we can apply Proposition \ref{correspondence} to conclude that the right cosets of $q_\ast(\pi_{1,q}^{\mathrm{coarse}}(X, b'))$ in $\pi_1^{\mathrm{coarse}}(X/G, b)$ are in bijection with $\sim$-equivalence classes of lifts of $b$. Note that we have a lift $g \circ b'$ of $b$ for every $g \in G$. It is easy to check that if $g \circ b' \sim h \circ b'$ then $g = h$ using the fact that the action is uniformly coarsely discontinuous. On the other hand, if $b''$ is some other lift of $b$ then for every $t$, $b''(t) = g_t \cdot b'(t)$ for some $g_t \in G$. If $|t - t'| \leq 1$, then $b'(t)$ and $g^{-1}_{t'} g_{t}(b'(t))$ are at most $2\rho(1)$ apart, where $\rho$ is the control function for $b$. But then $g_t = g_t'$ whenever $|t-t'| < 1$ outside of a bounded set. We have thus shown that for any lift $b''$ of $b$, we have $b'' \sim g \circ b'$ for a unique $g \in G$. This bijection between $G$ and $\sim$-classes of lifts of $b$ gives rise to a bijection from cosets of $q_\ast(\pi_{1,q}^{\mathrm{coarse}}(X, b'))$ in $\pi_1^{\mathrm{coarse}}(X/G, b)$ to elements of $G$ via Proposition \ref{correspondence}. All that remains is to check that this bijection respects the group operations and that the map $q_\ast$ is injective, both of which are easy to show.
\end{proof}

\begin{Corollary}\label{XGcoarsecor}
Let $G$ be a group acting on a metric space $X$ by isometries such that for every bounded set $K$, the set $\cup_{g \in G} g(K)$ is also bounded. Let $X/G$ be the orbit space and $q: X \to X/G$ the natural quotient map. Suppose further that the action of $G$ is uniformly coarsely discontinuous. Let $b: \mathbb{R}_+ \to X/G$ be an $\mathbb{R}_+$-basepoint such that for any two lifts $b', b'': \mathbb{R}_+ \to X$, there is a coarse homotopy $H$ from $b'$ to $b''$. Then for any lift $b'$ of $b$, we have a canonical short exact sequence
$$
\xymatrix{
0 \ar[r] & \pi_{1}^{\mathrm{coarse}}(X, b') \ar[r]^{q_\ast} &
 \pi_1^{\mathrm{coarse}}(X/G, b) \ar[r] & G \ar[r] &  0.
}
$$
\end{Corollary}
\begin{proof}
Easy once we notice that $q$ is coarse.
\end{proof}

Notice that the condition
$$K  \text{ bounded} \implies \cup_{g \in G} g(K) \text{ bounded}$$
from the above corollary is always satisfied if $G$ is finite.

Proposition \ref{correspondence}, Corollary \ref{correspondencecor}, Theorem \ref{XGcoarse1} and Corollary \ref{XGcoarsecor} all require the existence of coarse homotopies between lifts of $\mathbb{R}_+$-basepoints. This requirement is a natural analog of the requirement that the domain of a covering map be path-connected in order to get a classical lifting correspondence. In order to demonstrate that this condition is satisfied in a large number of examples, we connect it to a topological condition. The following lemma is based on Lemma 4.2 in \cite{HigsonRoeCoarseBC}, a result about maps from metric cones which used a different notion of coarse homotopy.

\begin{Lemma} \label{continuoustocoarse}
Let $X$ be a proper metric space and let $b, b': \mathbb{R}_+ \to X$ be two coarse and continuous maps. Suppose that there is a proper (classical) homotopy $h: \mathbb{R}_+ \times [0,1] \to X$ from $b$ to $b'$. Then there is a coarse homotopy from $b$ to $b'$. Moreover, if $f: X \to Y$ is a bornologous map such that $f \circ h$ is a proper map, then there is a $f$-coarse homotopy from $b$ to $b'$.
\end{Lemma}
\begin{proof}
It will be convenient to have $b(0) = b'(0)$ for the proof, so define new maps
$$
\overline{b} = \begin{cases}
h(0, 1/2-x) & x \leq 1/2 \\
b(x-1/2) & x \geq 1/2 
\end{cases}
$$
and 
$$
\overline{b}' = \begin{cases}
h(0, 1/2+x) & x \leq 1/2 \\
b'(x-1/2) & x \geq 1/2 
\end{cases}
$$
Note that $b$ and $b'$ are close to $\overline{b}$ and $\overline{b}'$ respectively (because $h(\{0\}\times [0,1])$ is compact and therefore bounded), so it is enough to show that $\overline{b}$ and $\overline{b}'$ are coarsely homotopic. Moreover, $\overline{b}$ and $\overline{b}'$ are continuously homotopic via a homotopy $\overline{h}$ with $\overline{h(0, t)} = \overline{b}(0) = \overline{b}'(0)$. We can transform $\overline{h}$ to be a proper continuous map from $c([0,1])$ as follows:
$$
\overline{h}_c(x,t) = \begin{cases}
\overline{h}(x, t/x) & x > 0 \\
\overline{h}(0,0) & x = 0
\end{cases}
$$
Since $\overline{h}_c$ is continuous, for every $K > 0$ there is an $L_K > 0$ so that
$$
d((x,t), (x',t')) \leq 1/L_K \implies d(f(x,t), f(x',t')) \leq 1
$$
for all $x, x', t, t' \in [0, K]$.  We may assume that the $L_K$ are increasing, are integer valued and are all at least $1$. We define a contraction $\rho: \mathbb{R}_+ \to \mathbb{R}_+$ piecewise as follows. Define $\rho$ on $[0,2L_2]$ to be the linear map from $[0, 2L_2]$ to $[0,1]$, and proceed by induction as follows: if $\rho$ has been defined for $x \in [0,n]$, $n \geq 1$, with $\rho(n) = k$, define $\rho$ on $[n, n+L_{k+2}(k+2)]$ to be an affine map with image $[k, k+1]$. Define $H(x,t) = \overline{h}_c(\rho(x), t\rho(x)/x)$. It is now easy to check that 
$$
d((x,t), (x',t')) \leq 1 \implies d(H(x,t), H(x',t')) \leq 1
$$
so that $H$ is coarse. It is not a coarse homotopy from $\overline{b}$ to $\overline{b}'$ because $H \circ i_0$ need not equal $\overline{b}$. However, the identity map on $\mathbb{R}_+$ and $\rho$ are coarsely homotopic by a straight-line homotopy, so the maps $\overline{b}$ and $H \circ i_0 = \overline{b}\circ \rho$ are coarsely homotopic. Similarly, so are $\overline{b}'$ and $H \circ i_1 = \overline{b}'\circ \rho$. Combining these homotopies with $H$ we obtain a coarse homotopy from $\overline{b}$ to $\overline{b}'$ as required. If $f: X \to Y$ is a bornologous map and $f\circ h$ is a proper map, then it is easy to check that $f$ composed with this coarse homotopy is also proper as constructed above.
\end{proof}

In the context of geometric group theory, a space $X$ (usually a CW-complex) is called \emph{semistable at $\infty$} if all proper continuous rays $f: \mathbb{R}_+ \to X$ are properly homotopic. It is still an open question whether all simply-connected spaces admitting a geometric group action by a finitely-presented one-ended group are semistable at $\infty$ (this is the so-called Semistability Conjecture; see for example \cite{Mihalik1983}). We refer the reader to Chapter 16 of \cite{Geoghegan} for a thorough treatment of semistability and its relation to the fundamental group at infinity. The following corollary shows that for length spaces, semistability implies ``coarse semistability'' (i.e.~that any two $\mathbb{R}_+$-basepoints are coarsely homotopic). 

\begin{Corollary} \label{semistability}
Let $X$ be a proper length space such that any two proper continuous rays $a, a': \mathbb{R}_+ \to X$ are properly homotopic, then any two coarse maps $b, b': \mathbb{R}_+ \to X$ are coarsely homotopic.
\end{Corollary}
\begin{proof}
Let $b, b': \mathbb{R}_+ \to X$ be two coarse maps. Define new maps $a$ and $a'$ as follows. For $x \in \mathbb{R}$, $a(x) = b(x)$ if $x$ is an integer, otherwise $a(x) = \gamma((x-n)\ell)$ where $\gamma$ is a path from $b(n)$ to $b(n+1)$ with length $\ell \leq d(b(n), b(n+1))+1$ and $n = \floor{x}$. Note that $\ell$ is bounded above since $b$ is coarse. Define $a'$ based on $b'$ in a similar fashion. The maps $a$ and $a'$ are continuous and are close to $b$ and $b'$ respectively. It follows that $a$ and $a'$ are both proper and that they are coarsely homotopic to $b$ and $b'$ respectively. Since $a$ and $a'$ are properly homotopic, they are also coarsely homotopic by Lemma \ref{continuoustocoarse}, so $b$ and $b'$ are coarsely homotopic as required.
\end{proof}

Even for spaces which are not semistable, their metric cones often are, as the following proposition shows. We recall the definition of metric or open cone on a manifold, which will be the central object of study in the next section.

\begin{Definition}[\cite{RoeWarped}]\label{conedef}
Given a Riemannian manifold $(M,g_M)$, the \textbf{metric cone} $\mathcal{O}M$ over $M$ is the manifold with boundary $M \times [1,\infty)$ with the Riemannian metric $t^2 g_M + g_\mathbb{R}$, where $t$ is the second coordinate and $g_M$ is the Riemannian metric on on $M$. We denote the geodesic metric on $M$ by $d_M$ and the geodesic metric on $\mathcal{O}M$ by $d_{\mathcal{O}M}$.
\end{Definition}

When $M$ is compact, this definition is known to coincide, up to coarse equivalence, with other common constructions of the metric cone (see e.g. \cite{RoeWarped}).

\begin{Proposition} \label{conesemistability}
Let $M$ be a Riemannian manifold and let $\mathcal{O}M$ be the open cone over $M$. Then any two coarse maps $b,b': \mathbb{R}_+ \to \mathcal{O}M$ are coarsely homotopic.
\end{Proposition}
\begin{proof}
Suppose for now that $b$ and $b'$ are also continuous maps. Since $M$ is connected, $\mathcal{O}M$ is clearly path-connected. Let $\alpha: [0,1] \to \mathcal{O}M$ be a path from $b(0)$ to $b'(0)$. Then there is a natural homotopy from $b$ to $b'$ given by 
$$
H(x,t) = \begin{cases}
b(x-3tx) & 0 \leq x \leq 1/3 \\
\alpha(3t-1) & 1/3 \leq x \leq 2/3 \\ 
b'(3tx-2x) & 2/3 \leq x \leq 1
\end{cases}
$$
This homotopy is far from coarse or even proper, so we make the following adjustment to obtain a new continuous proper homotopy $H'$:
$$
H'(x,t) = \left( p_1(H(x,t)),\ (1-t) \cdot p_2(H(x,0))+ t \cdot p_2(H(x,1)) \right)
$$
where $p_1: \mathcal{O}M \to M$ and $p_2: \mathcal{O}M \to \mathbb{R}_+$ are the projections to the two coordinates. By Lemma \ref{continuoustocoarse}, this proper continuous homotopy induces a coarse homotopy from $b$ to $b'$. The result now follows from Corollary \ref{semistability}.
\end{proof}

We can now exhibit a class of examples satisfying the requirement that lifts of $\mathbb{R}_+$-basepoints are coarsely homotopic.

\begin{Corollary}\label{coneslifting}
Let $M$ and $N$ be connected Riemannian manifolds with $N$ compact. Let $f: \mathcal{O}M \to \mathcal{O}N$ be a map on metric cones of the form $f(x,t) = (f_1(x,t),t)$. Then any two $f$-coarse maps $b,b': \mathbb{R}_+ \to \mathcal{O}M$ are $f$-coarsely homotopic.
\end{Corollary}
\begin{proof}
Notice that the coarse homotopy constructed in Proposition \ref{conesemistability} is $p_2$-coarse, where $p_2$ is the project of $\mathcal{O}M$ onto the second coordinate. Since $N$ is compact, and thus bounded, it follows that this coarse homotopy is also $f$-coarse.
\end{proof}

In the next section, we will explore maps between metric cones further to show how covering maps between manifolds give rise to soft quotient maps with scattered fibres between their metric cones. This will allow us to compute the coarse fundamental groups of metric cones over compact, non-positively curved manifolds.

\section{Application to metric cones}\label{cones}
In order to demonstrate how the results of the previous section can be used to compute coarse fundamental groups of spaces, we consider metric cones over compact Riemannian manifolds of non-positive sectional curvature (for convenience, all manifolds in this section will be assumed to be connected unless stated otherwise). Metric cones were defined at the end of the previous section (Definition \ref{conedef}). Specifically, we prove the following theorem. 

\begin{Theorem}\label{main}
Let $M$ be a compact Riemannian manifold of non-positive sectional curvature, and let $\mathcal{O}M$ be the metric cone over $M$. Let $b$ be any $\mathbb{R}_+$-basepoint in $\mathcal{O}M$. Then there is an isomorphism
$$
\pi_{1}^\mathsf{coarse}(\mathcal{O}M, b) \cong \pi_1(M),
$$
where $\pi_1(M)$ is the fundamental group of $M$.
\end{Theorem}

The main idea of the proof is to write such a cone as the quotient of the cone over the universal cover of the manifold. This works because the cone over the universal cover is sufficiently ``trivial'' at the level of coarse fundamental groups. Unfortunately, we require a number of technical lemmas before we can execute this idea. We first need some basic results on distances in $\mathcal{O}M$. 

\begin{Lemma}\label{conetech}
Let $\mathcal{O}M$ be the metric cone over a complete Riemannian manifold $M$. If $(x,t)$ and $(x',t')$ are two points in $\mathcal{O}M$, then we have the following inequalities.
\begin{align*}
\mathrm{(1)} \quad & d_{\mathcal{O}M}((x,t), (x',t')) \leq |t-t'| + d_M(x,x')\cdot \min(t,t') \\
\mathrm{(2)} \quad & |t-t'| \leq  d_{\mathcal{O}M}((x,t), (x',t')) \\
\mathrm{(3)} \quad & d_M(x,x')  \leq  d_{\mathcal{O}M}((x,t), (x',t')) \\
\mathrm{(4)} \quad & \forall R\ \exists S \ d_{\mathcal{O}M}((x,t), (x',t')) < R \implies d_M(x,x')\cdot t \leq S
\end{align*}
\end{Lemma}
\begin{proof}
\begin{enumerate}
\item[(1)] Suppose without loss of generality that $t \leq t'$ and consider the concatenation of the path $\alpha: [0, t'-t] \to \mathcal{O}M$ given by $\alpha(s) = (x',t'-s)$ and $\beta(s) = (\gamma(s), t)$ where $\gamma$ is the geodesic from $x'$ to $x$ in $M$. This is a unit-speed path with length $|t-t'| + d_M(x,x')t$ from which we obtain the result.
\item[(2)] The projection $p_2: \mathcal{O}M \to \mathbb{R}_+$ onto the second coordinate is distance non-increasing because the metric on $\mathcal{O}M$ is $t^2g_M+g_\mathbb{R}$. The result follows. 
\item[(3)] The projection $p_1: \mathcal{O}M \to M$ onto the first coordinate is also distance non-increasing because the metric on $\mathcal{O}M$ is $t^2g_M+g_\mathbb{R}$ and $t \geq 1$. The result follows. 
\item[(4)] Let $\varepsilon > 0$. Let $R > 0$ and let $d_{\mathcal{O}M}((x,t), (x',t')) < R$. Let $\gamma$ be a piecewise differentiable curve in $\mathcal{O}M$ from $(x,t)$ to $(x',t')$ with length at most $R$. By inequality (2), the image of $\gamma$ is contained in $[\max(1,t-R), t+R] \times M$. Choosing a decomposition of the tangent space into an $M$ and an $\mathbb{R}$ subspace, we have
\begin{align*}
R  \ & \geq \int ||\gamma'(s)|| ds \\
& \geq \int ||\gamma'_M(s)|| ds \\
& \geq  \int ||\min(1,t-R) (p_1 \circ \gamma)' (s)|| ds \\
& \geq \min(1,t-R) d_M(x,x')
\end{align*}
where $\gamma'_M(s)$ is the $M$-component of the tangent vector $\gamma'(s)$. This gives 
$$d_M(x,x') \cdot t \leq \frac{Rt}{\max(1,t-R)}$$
which is bounded above independent of $t$. This completes the proof.
\end{enumerate}
\end{proof}

The main fact about Riemannian manifolds of non-positive sectional curvature we will need is the following lemma. Recall that a \textbf{Hadamard manifold} is a complete simply connected Riemannian manifold of non-positive sectional curvature.

\begin{Lemma}\label{lip}
Let $\gamma_1, \gamma_2: [0, \infty) \to M$ be unit speed geodesics on a Hadamard manifold. Then for any $\theta \in [0,1]$ and $q,q' > 0$,
$$d(\gamma_1(\theta q), \gamma_2(\theta q')) \leq \max(d(\gamma_1(0), \gamma_2(0)), d(\gamma_1(q), \gamma_2(q'))).$$
\end{Lemma}
\begin{proof}
This follows from the fact that any Hadamard manifold is a CAT(0) space (this is a consequence of the Cartan-Hadamard theorem and a result of Alexandrov;  see Theorem 1A.6 and Theorem 4.1 in \cite{BridsonHaefliger}). 
\end{proof}

The following lemmas allow us to move from locally Lipschitz to globally Lipschitz via coarse homotopy.

\begin{Lemma}[Lemma 1.9 in \cite{MitchenerNorSchick}]\label{pastinglip}
Let $X$ be a geodesic space, and let $X = \cup_i A_i$ be a decomposition of $X$ into closed subsets so that every compact set in $X$ intersects only finitely many of the $A_i$. If $f: X \to Y$ is a map which is $C$-Lipschitz when restricted to each $A_i$, then $f$ is $C$-Lipschitz.
\end{Lemma}

\begin{Lemma}\label{localtoglobal}
Let $f: c([0,1]) \to Y$ be a map to a metric space $Y$ such that for every $K > 0$ there is an $L_K > 0$ so that
$$
d(f(x,t), f(x',t')) \leq L_K\cdot d((x,t), (x',t'))
$$
for all $x, x', t, t' \in [0, K]$. Then there is a coarse map $g: c([0,1]) \to c([0,1])$ which is coarsely homotopic to the identity such that $f \circ g$ is a $1$-Lipschitz map. In particular, $f$ is coarsely homotopic to a $1$-Lipschitz map.
\end{Lemma}
\begin{proof}
We may assume that the $L_K$ are increasing, are integer valued and are all at least $1$. The map $g$ will have the form $g(x,t) = (\rho(x), t\rho(x)/x)$ for a monotone map $\rho: \mathbb{R}_+ \to \mathbb{R}_+$. Define $\rho$ on $[0,2L_2]$ to be the linear map from $[0, 2L_2]$ to $[0,1]$, and proceed by induction as follows: if $\rho$ has been defined for $x \in [0,n]$, $n \geq 1$, with $\rho(n) = k$, define $\rho$ on $[n, n+L_{k+2}(k+2)]$ to be an affine map with image $[k, k+1]$. One easily checks that the map $f \circ g$ is $1$-Lipschitz using Lemma \ref{pastinglip}. Moreover, $g$ is coarsely homotopic to the identity via the obvious straight line homotopy.
\end{proof}

\begin{Lemma}\label{coarsetolip}
Let $f: c([0,1]) \to \mathcal{O}M$ be a $p_2$-coarse map where $M$ is a Hadamard manifold with $p_2:\mathcal{O}M \to [1, \infty)$ the projection onto the second coordinate. Then $f$ is $p_2$-coarsely homotopic to a $1$-Lipschitz map $f'$. Moreover, if $f \circ i_0 = f \circ i_1$, then we can choose the $p_2$-coarse homotopy $H$ so that $H((x,0),s) = H((x,x),s)$ for all $s$.
\end{Lemma}
\begin{proof}
Let $A = c([0,1]) \cap \mathbb{Z}^2$ be the set of points in $c([0,1])$ with integer coordinates, and let $f_A$ be the restriction of $f$ to $A$. Each of the maps $p_1 \circ f_A: A \to M$ and $p_2 \circ f_A: A \to [1, \infty)$, where $p_1$ is the projection onto the first coordinate, can be extended to Lipschitz maps by interpolating with geodesics. Taking the product of these extensions, we obtain a map $g: c([0,1]) \to \mathcal{O}M$ which using Lemma \ref{conetech} can be shown to be close to $f$ (and hence coarsely homotopic to it), but which may be only locally Lipschitz on each square $c([0,1]) \cap [k, k+1]^2$. However, by Lemma \ref{localtoglobal}, $g$ is $p_2$-coarsely homotopic to a $1$-Lipschitz map $f'$ on all of $c([0,1])$, so we obtain the required result (after checking that both coarse homotopies satisfy the additional condition at the end of the lemma).
\end{proof}

The proof of the following Proposition is based on the proof of Proposition 5.2 in \cite{MitchenerNorSchick}. 

\begin{Proposition}\label{contraction}
Let $\alpha: c([0,1]) \to \mathcal{O}M$ be a $p_2$-coarse map where $M$ is a Hadamard manifold with $p_2$ the projection onto the second coordinate. Suppose that $\alpha \circ i_0 = \alpha \circ i_1$. Then $\alpha$ is $p_2$-coarsely homotopic to a map of the form $\beta(x,t) = b(x)$ via a $p_2$-coarse homotopy $H$ satisfying $H((x,0),s) = H((x,x),s)$ for all $s$.
\end{Proposition}
\begin{proof}
Pick any $p \in M$. By adjusting $\alpha$ at one point if necessary, we may assume that $\alpha(0,0) = (p,1)$. By Lemma \ref{coarsetolip}, we may also suppose that $\alpha$ is $1$-Lipschitz. We start by noticing that the map $\alpha$ is coarsely homotopic to the map $\alpha': c([0,1]) \to \mathcal{O}M$ given by 
$$\alpha'(x,t) = \alpha(x/\max(1,\sqrt{x}), t/\max(1,\sqrt{x})).$$
via a homotopy $H_1((x,t),s) = (x-s, t(x-s)/x)$ for $s \leq x- \max(1,\sqrt{x}))$ composed with $\alpha$. Note that $\alpha \circ H_1((x,0),s) = \alpha \circ H_1(\alpha(x,x), s)$ for all $s$. Now consider the map $H: I_q \mathcal{O}M \to \mathcal{O}M$ defined by $H((y,u),s) = (f(y,s/u), u)$, with $f$ is defined by
$$
f(x,t) = \begin{cases}
\gamma_x(t) & t \leq d(p,x) \\
p & t > d(p,x) \\
\end{cases}
$$
where $\gamma_x$ is the unique unit speed geodesic from $x$ to $p$ in $M$, and $q(x,t) = d(x,p)t-1$. This naturally leads to a map
$$
H': I_{ q \circ \alpha'} c([0,1]) \to \mathcal{O} M
$$
given by $H'((x,t), s) = H(\alpha'(x, t), s)$ for $s \leq q(\alpha'(x,t)))$. Note that $H' \circ i_1$ has image contained in $p \times [1,\infty)$. We now show that $H'$ is a $p_2$-coarse homotopy.

Let $(x,t)$ and $(x',t')$ be two points of distance at most $1$ apart in $c([0,1])$, with $x,x' > 1$. Using the Mean Value Theorem, we have
\begin{align*}
d((\sqrt{x}, t/\sqrt{x}), (\sqrt{x'}, t'/\sqrt{x'})) & \leq |\sqrt{x} - \sqrt{x'}| + |t/\sqrt{x} - t/\sqrt{x'}| + |t/\sqrt{x'} - t'/\sqrt{x'}| \\
& \leq \frac{1}{2\sqrt{x}} + \frac{t}{2x\sqrt{x}} + \frac{|t-t'|}{\sqrt{x'}}  \leq \frac{2}{\sqrt{x}} 
\end{align*}
Let $\alpha'(x,t) = (y,u)$ and $\alpha'(x',t') = (y',u')$. By the above and Lemma \ref{conetech}, we have $|u-u'| \leq 2/\sqrt{x}$ and $d(y,y') \leq 2/\sqrt{x}$. Since $\alpha'(0,0) = (p,1)$, we have that $d(y,p)$ and $u$ are both bounded above by $(t/\sqrt{x} + \sqrt{x}) \leq 2\sqrt{x}$ by Lemma \ref{conetech}. Thus,
\begin{align*}
|q(y,u) - q(y',u')| &  \leq d(y',p)\cdot |u'-u| + u\cdot |d(y,p) - d(y',p)| \\
& \leq d(y',p)\cdot |u'-u| + u\cdot d(y,y') \\
& \leq 2\sqrt{x}\cdot \frac{2}{\sqrt{x}} + 2\sqrt{x}\cdot \frac{2}{\sqrt{x}} =8 
\end{align*}
which, using the fact that $c([0,1])$ is a geodesic space, shows that $q\circ \alpha'$ is coarse for the region where $x,x' > 1$. Since the image of the region where $x,x' \leq 1$ under $q\circ \alpha'$ is clearly bounded, we get coarseness on all of $c([0,1])$. With a view to showing coarseness of $H'$, suppose further that $s \leq \min(q(\alpha'(x,t))), q(\alpha'(x',t')))$. The distance between $H'((x,t),s)$ and  $H'((x',t'), s))$ is bounded above by 
$$
d(H((y,u),s), H((y',u), s)) + d(H((y',u),s), H((y',u'), s)).
$$
Looking at the triangle $(p,y,y')$ and using the fact that $M$ is a CAT(0) space, we have that $d(f(y, s/u), f(y', s/u)) \leq d(y,y')$, so the first term is bounded by
$$
u \cdot d(y,y')\leq (2 \sqrt{x}) (2/ \sqrt{x}) = 4.
$$ 
For the second term, assume without loss of generality that $u \leq u'$. Then the second term is bounded above by (using the Mean Value theorem again),
\begin{align*}
u \cdot d(f(y',s/u), f(y', s/u')) + |u-u'| & \leq u\cdot s \cdot (1/u - 1/u') + |u-u'| \\
& \leq u \cdot u \cdot d(y,p) \cdot \frac{2}{\sqrt{x}}\cdot \frac{1}{u^2} + \frac{2}{\sqrt{x}}\\
& \leq d(y,p)\cdot \frac{2}{\sqrt{x}} + \frac{2}{\sqrt{x}} \leq 4 + \frac{2}{\sqrt{x}}\\
\end{align*}
which is bounded. It is easy to check that $H'$ is coarse in the third coordinate; indeed,
$$
d(H((y,u),s), H((y,u), s')) \leq |s/u-s'/u|*u \leq |s-s'|.
$$
We have thus shown that $H'$ is coarse, and it is moreover clearly $p_2$-coarse if $\alpha$ is $p_2$-coarse. The image of the map $\beta' = H' \circ i_1$ is completely contained in $p \times [1,\infty)$, and so is $p_2$-coarsely homotopic relative to the boundary $\partial c([0,1])$ to the map $\beta(x,t) = \beta'(x,0)$ by a linear homotopy (or by invoking Theorem 5.6 of \cite{MitchenerNorSchick} for the space $[1,\infty)$). It is easy to check that $H'((x,0),s) = H'((x,x),s)$ for all $s$, from which the last condition in the statement follows.
\end{proof}

We are almost done with technical proofs; the following result brings together the previous lemmas to show that that the metric cone over $\tilde{M}$ is trivial at the level of coarse homotopy.

\begin{Corollary}\label{trivialM}
Let $M$ be a compact Riemannian manifold of non-positive sectional curvature and let $\tilde{M}$ be its universal cover equipped with the metric lifted from $M$. Let $\sigma: \mathcal{O}\tilde{M} \to \mathcal{O}M$ be the map on cones induced by the covering map $\tilde{M} \to M$. Then 
\begin{itemize}
\item[(1)] any two $\sigma$-coarse $\mathbb{R}_+$-basepoints $b$ and $b'$ are $\sigma$-coarsely homotopic.
\item[(2)] for any $\sigma$-coarse $\mathbb{R}_+$-basepoint $b$ in $\mathcal{O}\tilde{M}$, $\pi_{1,\sigma}^\mathsf{coarse}(\mathcal{O}\tilde{M}, b)$ is trivial. 
\end{itemize}
\end{Corollary}
\begin{proof}
Statement (1) follows from Corollary \ref{coneslifting}. Let $\alpha$ represent a class in $\pi_{1,\sigma}^\mathsf{coarse}(\mathcal{O}\tilde{M}, b)$. We may replace $\alpha$ up to closeness (and hence up to relative $p_2$-coarse homotopy) by a map which is constant in $t$ on $[0,1)^2 \cap c([0,1])$. Proposition \ref{contraction} shows that $\alpha$ is $\sigma$-coarsely homotopic via $H$ to a ``trivial'' map $H\circ i_1(x,t) = b'(x)$. Again, make the following small adjustment for technical reasons which does not affect the coarseness of $H$: redefine $H$ so that $H((x,t), s) = H((x,t), 0)$ for $(x,t) \in [0,1)^2$. Having done this, we can use an easy adaptation of Lemma 2.6 in \cite{MitchenerNorSchick} to assume that $H$ is of the form $H: I_p c([0,1]) \to \mathcal{O}M$ with $p(x,t) = x-1$. The coarse homotopy $H$ is still not a homotopy relative to the boundary, though, so we have to use a trick which is familiar from topology. Construct a new map $H': I_{p(x)} c([0,1]) \to \mathcal{O}M$ as follows:
$$
H'((x,t),s) = 
\begin{cases}
b(x) & t \leq x/4 - s/4 \\
H((x,0),4t-x+s) & x/4-s/4 \leq t \leq x/4 \\
 H((x,2t-x/2),s) & x/4  \leq t \leq 3x/4    \\
 H((x,0),3x-4t+s) & 3x/4 \leq t \leq 3x/4+s/4 \\
 b(x) & 3x/4 +s/4 \leq t \leq x \\
\end{cases}
$$
It is straightforward to show that this is a $p_2$-coarse homotopy relative to the boundary from $\alpha$ to a concatenation of a map $\lambda$ from $c([0,1])$, followed by map constant in $t$, followed by the reverse of $\lambda$. This concatenation is easily shown to be $p_2$-coarsely homotopic to the trivial element of $\pi_{1,\sigma}^\mathsf{coarse}(\mathcal{O}\tilde{M}, b)$, which completes the proof of (2).
\end{proof}

\begin{Lemma}\label{discon}
Let $M$ be a Riemannian manifold and let $G$ be a group acting on $M$ properly discontinuously and cocompactly by isometries. Then the induced action of $G$ on $\mathcal{O}M$ is uniformly coarsely discontinuous.
\end{Lemma}
\begin{proof}
Since $G$ acts properly discontinuously by isometries, there is a global $C > 0$ so that $d(x, g\cdot x) \geq C$ for all $x \in M$ and all $g \in G \setminus \{e\}$. Let $R > 0$ and let $g \in G \setminus \{e\}$. By (4)  in Lemma \ref{conetech}, there exists a $T$ such that if $t > T$ then the distance between $(x,t)$ and $(g\cdot x, t)$ is at least $R$. Since the action is cocompact, $\mathcal{O}M \cap M \times [0, T]$ is contained in $\cup_{g\in G} g(K)$ for some compact $K$, which proves the lemma.
\end{proof}

\begin{proof} [Proof of Theorem \ref{main}]
Consider the map $\sigma: \mathcal{O}\tilde{M} \to \mathcal{O}M$ induced by the covering map from the universal cover $\tilde{M}$. Since the metric on $\mathcal{O}\tilde{M}$ is lifted from $\mathcal{O}M$, $\mathcal{O}M$ is isometric to the quotient of $\mathcal{O}\tilde{M}$ by the action of $\pi_1(M)$ induced by the action of $\pi_1(M)$ on $\tilde{M}$ (with the metric given as in Example \ref{exgroup}). By Lemma \ref{discon}, the action of $\pi_1(M)$ is uniformly coarsely discontinuous, so the result follows from Theorem \ref{XGcoarse1}  and Corollary \ref{trivialM}.
\end{proof}

In Theorem 5.6 of \cite{MitchenerNorSchick}, MNS prove a similar result for cones over finite simplicial complexes. This suggests that the curvature condition could possibly be relaxed in Theorem \ref{main} above. This question, as well as the question of whether Theorem \ref{main} can be recovered from the result in \cite{MitchenerNorSchick} via a triangulation argument, is left for a future paper. Even if Theorem \ref{main} is a corollary to the result for simplicial complexes, it still serves as an illustrative example of computing coarse fundamental groups using the Coarse Lifting Lemma


\begin{thebibliography}{99}










\bibitem{BridsonHaefliger} M. R. Bridson, and A. Haefliger, Metric spaces of non-positive curvature, Vol. 319, Springer Science \& Business Media, 2013.



\bibitem{BunkeEngel} U. Bunke and A. Engel, \emph{Homotopy theory with bornological coarse spaces}, preprint, 	arXiv:1607.03657. 



\bibitem{Dranishnikov00} A. N. Dranishnikov, \emph{Asymptotic topology}, Russian Mathematical Surveys 55.6, 2000, 1085.




\bibitem{DikranjanZava} D. Dikranjan and N. Zava, \emph{Some categorical aspects of coarse spaces and balleans}, Topology and its Applications 225, 2017, 164--194.



\bibitem{DH} J. Dydak and C.S. Hoffland, \emph{An Alternative Definition of Coarse Structures}, Topology and its Applications 155 (9), 2008, 1013-1021.






\bibitem{Geoghegan} R. Geoghegan, \emph{Topological methods in group theory}, Springer Science \& Business Media, 2007.


\bibitem{Gromov93} M. Gromov, Asymptotic invariants of infinite groups, Geometric group theory vol. 2 (Sussex, 1991), London Math. Soc. Lecture Note Ser. vol. 182, Cambridge University Press, Cambridge, 1993, 1--295.



\bibitem{WeighillHigginbotham} L. Higginbotham and T. Weighill, \emph{Coarse quotients by group actions and the maximal Roe algebra}, Journal of Topology and Analysis (online ready), https://doi.org/10.1142/S1793525319500341.

\bibitem{HigsonRoeCoarseBC} N. Higson and J. Roe, \emph{On the coarse Baum-Connes conjecture}, Novikov conjectures, index theorems and rigidity 2, 1995, 227--254.

\bibitem{Mihalik1983} Mihalik, M. L, \emph{Semistability at the end of a group extension}, Transactions of the American Mathematical Society 277.1, 1983, 307--321.









\bibitem{MitchenerNorSchick} P. Mitchener, B. Norouzizadeh and T. Schick, \emph{Coarse homotopy groups}, arXiv preprint, arXiv:1811.10096.

\bibitem{Munkres} J. Munkres, {\em Topology}, Prentice Hall, 2000, 2nd edition.





\bibitem{RoeComplete} J. Roe, Coarse cohomology and index theory on complete Riemannian manifolds, Memoirs of the American Mathematical Society 497, 1993.

\bibitem{RoeIndex} J. Roe, Index theory, coarse geometry, and topology of manifolds, CBMS Regional Conference Series in Mathematics Volume 90, Published for the Conference Board of the Mathematical Sciences, Washington, DC, 1996.


\bibitem{Roe} J. Roe, {\em Lectures in Coarse Geometry}, University Lecture Series 31, American Mathematical Society, Providence, RI, 2003. 

\bibitem{RoeWarped} J. Roe, \emph{Warped cones and property A}, Geometry \& Topology 9.1, 2005, 163--178.





\bibitem{WillettYu} R. Willett and G. Yu, \emph{Higher index theory for certain expanders and Gromov monster groups, I}, Advances in Mathematics 229, 2012, 1380--1416. 



 



\bibitem{YuCoarseBC} G. Yu, \emph{The coarse Baum-Connes conjecture}, K-theory 9, 1995.

\bibitem{YuNovikov} G. Yu, \emph{The Novikov conjecture for groups with finite asymptotic dimension}, Annals of Mathematics 147, 1998, 325--355.


\bibitem{YuEmbed} G. Yu, \emph{The coarse Baum-Connes conjecture for spaces which admit a uniform embedding into Hilbert space}, Inventiones mathematicae 139.1, 2000, 201-240.


\end{thebibliography}
\end{document}